# Composition Direction of Seymour's Theorem for Regular Matroids — Formally Verified

A Preprint

**Martin Dvorak**
ISTA
Klosterneuburg, Austria
martin.dvorak@matfyz.cz

**Tristan Figueroa-Reid**
Reed College
Portland, United States of America
pub.tristanf@gmail.com

**Rida Hamadani**
LMAP, UPPA
Pau, France
mridahamadani@gmail.com

**Byung-Hak Hwang**
Korea Institute for Advanced Study
Seoul, South Korea
byunghakhwang@gmail.com

**Evgenia Karunus**
University Of Bonn
Bonn, Germany
lakesare@gmail.com

**Vladimir Kolmogorov**
ISTA
Klosterneuburg, Austria
vnk@ist.ac.at

**Alexander Meiburg**
University of Waterloo & Perimeter Institute of Theoretical Physics
Waterloo, Canada
srmdt@ohaithe.re

**Alexander Nelson**
thmprover@gmail.com

**Peter Nelson**
University of Waterloo
Waterloo, Canada
apnelson@uwaterloo.ca

**Mark Sandey**
UC Riverside
Riverside, United States of America
mark@sandey-family.com

**Ivan Sergeev**
ISTA
Klosterneuburg, Austria
i.i.sergeyev@gmail.com

## Abstract

Seymour's decomposition theorem is a hallmark result in matroid theory presenting a structural characterization of the class of regular matroids. Formalization of matroid theory faces many challenges, most importantly that only a limited number of notions and results have been implemented so far. In this work, we formalize the proof of the forward (composition) direction of Seymour's theorem for regular matroids. To this end, we develop a library in Lean 4 that implements definitions and results about totally unimodular matrices, vector matroids, their standard representations, regular matroids, and 1-, 2-, and 3-sums of matrices and binary matroids given by their standard representations. Using this framework, we formally state Seymours decomposition theorem and implement a formally verified proof of the composition direction in the setting where the matroids have finite rank and may have infinite ground sets.

*We would like to dedicate the paper to the memory of Klaus Truemper, whose monograph Matroid Decomposition [12] laid the foundation for our entire work.*

## 1 Introduction

Seymour's regular matroid decomposition theorem is a hallmark structural result in matroid theory [9, 12, 4, 7]. It states that, on the one hand, any 1-, 2-, and 3-sum of two regular matroids is regular, and on the other hand, any regular matroid can be decomposed into matroids that are graphic, cographic, or isomorphic to $R_{10}$ by repeated 1-, 2-, and 3-sum decompositions.

The interest in matroids comes from the fact that they capture and generalize many mathematical structures and properties, such as linear independence (captured by vector matroids), graphs (graphic matroids), and extensions of fields (algebraic matroids). Another advantage of matroids is that they admit a relatively short definition, making them amenable to formalization. As for Seymour's theorem, it not only presents a structural characterization of the class of regular matroids, but also leads to several important applications, such as polynomial algorithms for testing if a matroid is binary and for testing if a matrix is totally unimodular. Additionally, Seymour's theorem can offer a structural approach for solving certain combinatorial optimization problems, for example, it leads to the characterization and efficient algorithms for the cycle polytope.

Formalization of results about matroids faces several challenges. One of them is that the support for them is limited. In Mathlib, only selected basic definitions for matroids are implemented, such as maps, duals, and minors. However, many other fundamental notions are not yet implemented, including representability and regularity, the splitter theorem and the separation algorithm. Part of the difficulty stems from the fact that classically, matroids are defined only in the finite case (i.e., when the ground set and the rank are finite), while Mathlib implements matroids more generally, allowing them to be infinite and to have infinite rank. Additionally, the proofs presented in the existing literature require substantial additional work to make them easily amenable to formalization.

The goal of our work was to develop a general and reusable library proving a result that is at least as strong as the forward (composition) direction of classical Seymour's theorem (i.e., stated for finite matroids). Moreover, our aim was to make our library modular and extensible by ensuring compatibility with matroids in Mathlib [8].

To achieve our goals, we made the following compromises. First, we focused on the implementation of the proof of the composition direction, while only stating the decomposition direction. Second, we assumed finiteness where it would simplify proofs, while making sure that the final results held for finite matroids (in fact, they hold for matroids with potentially infinite ground set and finite rank). Finally, we tailored our implementation specifically to Seymour's theorem, avoiding introducing additional matroid notions if possible. Our project makes the following contributions:

- Formalized definition and selected properties of totally unimodular matrices, some of which were added to Mathlib.
- Implemented definitions and formally proved selected results about vector matroids, their standard representations, regular matroids, and 1-, 2-, and 3-sums of matrices and vector matroids given by their standard representations.
- Implemented a formally verified proof of the composition direction of Seymour's theorem, i.e., that any 1-, 2-, and 3-sum of two regular matroids is regular, in the case where the matroids may have infinite ground sets and have finite rank.
- Stated the decomposition direction of Seymour's theorem, i.e., that any regular matroid of finite rank can be decomposed into matroids that are graphic, cographic, or isomorphic to $R_{10}$ by repeated 1-, 2-, and 3-sum decompositions.

Our formalization[1] is conceptually split into two parts: "implementation" and "presentation". Implementation is contained in the `Seymour` folder and encompasses all definitions and lemmas used to obtain our results. Presentation is contained in the `Seymour.lean` file, which repeats selected definitions and theorems comprising the key final results of our contribution. Every definition in the "presentation" file is checked to be definitionally equal to its counterpart from the "implementation" using the `recall` or the `example` command. Similarly, we `recall` every theorem presented here and then use the `#guard_msgs in #print axioms` command to check that the implementation of its proof (including the entire dependency tree) depends only on the three axioms `[propext, Classical.choice, Quot.sound]`, which are standard for Lean projects that use classical logic.

We refer to the statements of the final results and the definitions they (transitively) depend on as *trusted code*. The `Seymour.lean` file repeats all nontrivial trusted code, so that the reader can believe [10] our results without having to examine the entire implementation, assuming that the reader also uses the Lean compiler to check that all proofs are correct. Note that basic definitions from Lean and Mathlib are part of the trusted code but are not repeated in `Seymour.lean`, and we let the reader decide whether to blindly trust them or read them as well.

[1] `https://github.com/Ivan-Sergeyev/seymour/tree/v1`

While working on our project, we leveraged the LeanBlueprint[2] tool to help guide our formalization efforts. In particular, we used it to create theoretical blueprints and dependency graphs, which allowed us to get a clearer overview of the results we were formalizing, as well as their dependencies. In our workflow, we first created a write-up encompassing the classical results from [12]. Based on this write-up, we developed a self-contained theoretical blueprint for our formalization by filling in gaps, fleshing out technical details, and sometimes re-working certain proofs. We followed this blueprint during the development of our library, keeping it up to date and turning it into documentation of our code.

We use Lean version 4.18.0 and we import Mathlib library revision aa936c3 (dated 2025-04-01).

We made the code snippets in this paper as faithful to the content of the repository as possible, though we made some omissions. In particular, proofs inside definitions were replaced by the `sorry` keyword in the paper, while the repository contains full implementation.

## 2 Theory Underpinning the Formalization

There are two classical sources presenting the proof of Seymour's decomposition theorem: [9] and [12], each with their own advantages and disadvantages.

Oxley [9] develops a general theory of matroids and has a broader focus. It introduces many abstract notions and proves many statements about them, and Seymour's theorem and its dependencies are also stated and proved in terms of these abstract notions alongside many other results. The advantages of following [9] would be the higher reusability, generality, and extensibility of the formalization. Indeed, since [9] introduces a lot of foundational notions and results, the resulting implementation could serve as the basis for formalization of many other results from classical matroid theory. Moreover, [9] is more general than [12] in certain aspects, for example, [12] defines 1- and 2-sums only for binary matroids, while their definitions in [9] do not have this restriction. Finally, it seems that the approach to theory of infinite matroids [3] is more closely aligned with the approach of [9] than [12], which might make it easier to generalize formalizations based on the former than the latter to the infinite matroid setting. However, proof formalization following [9] would face many challenges. First, the support for matroids in Mathlib [8] at the time we carried out our project was quite limited. Thus a lot of time would be dedicated to developing low-level definitions and results about them, especially in the infinite matroid setting to ensure compatibility with Mathlib. Second, certain intermediate results could turn out difficult to formally prove. From our experiments, proving the equivalence of multiple characterizations of regular matroids turned out hard to formalize. Finally, [9] leaves many technical steps as exercises for the reader, most crucially leaving out the proof of regularity of 3-sum, and contains many proofs that crucially rely on graph theory which was not supported in Mathlib. This would make it challenging to convert the proofs to their formalized versions.

In contrast, Truemper [12] focuses on decomposition and composition of matroids, with Seymour's theorem being one of the most prominent theorems that it builds towards. Truemper [12] more frequently than Oxley [9] utilizes explicit matrix representations in definitions, theorems, and proofs, especially when it comes to 1-, 2-, and 3-sums of regular matroids. Thus, following [12] would require implementing fewer intermediate definitions and results to begin working with Seymour's theorem itself. Moreover, Mathlib's support for matrices and linear independence was more extensive than for matroids, so this would allow us to build upon more things that were already available. However, following the approach of [12] had several important limitations. As mentioned earlier, it would be less general and potentially less amenable to generalization to the infinite matroid setting than [9]. Moreover, faithfully following [12] would mean implementing similar definitions and theorems on several levels of abstraction. More specifically, 1-, 2-, and 3-sums would need to be implemented separately for matrices, binary matroids defined by standard representation matrices, and binary matroids in general, and the results about the sums of these objects would need to be proved and propagated accordingly. Last but not least, similar to [9], one would need to fill in the omitted technical details and re-work proofs that could be extremely challenging to formalize directly, especially those involving graph-theoretic arguments.

Ultimately, we decided to follow the approach of [12] over [9] for formalizing Seymour's theorem, as it aligned more closely with our goals and values. We aimed to formalize the statement of Seymour's theorem and the proof of the composition direction, so having to implement fewer intermediate definitions and lemmas and being able to use more tools from Mathlib was a big plus. Though we did not mind limiting the generality of our contributions to classical results, our final results go beyond that and hold for matroids of finite rank with potentially infinite ground sets. The completeness of the presentation in [12] allowed us to develop a theoretical blueprint, where we fleshed out the technical details, circumvented problematic intermediate results, and streamlined the proofs, especially in the case of 3-sums.

[2] `https://github.com/PatrickMassot/leanblueprint`

## 3 Preliminaries

Throughout this paper, we write $\mathbb{Z}_n$ to denote `ZMod n` for any positive integer $n$, most often in the case $\mathbb{Z}_2$ denoting `ZMod 2`, which is also written as `Z2` in the code.

This section reviews Mathlib declarations our code relies on.

### 3.1 Matroids

Matroids have many equivalent definitions [9, 12, 3]. In Mathlib, the structure `Matroid` captures the definition via the *base conditions* from [3]: a *matroid* is a pair $M = (E, \mathcal{B})$ where $E$ is a (potentially infinite) ground set and $\mathcal{B} \subseteq 2^E$ is a collection of sets such that:

(i) $\mathcal{B} \neq \emptyset$.

(ii) For all $B_1, B_2 \in \mathcal{B}$ and all $b_1 \in B_1 \setminus B_2$, there exists $b_2 \in B_2 \setminus B_1$ such that $(B_1 \setminus \{b_1\}) \cup \{b_2\} \in \mathcal{B}$.

(iii) For all $X \subseteq E$ and $I \subseteq X$ such that $I \subseteq B_1$ for some $B_1 \in \mathcal{B}$, there exists a maximal $J$ such that $I \subseteq J \subseteq X$ and $J \subseteq B_2$ for some $B_2 \in \mathcal{B}$.

A set $B \in \mathcal{B}$ is called a *base*, and (ii) is known as the *base exchange property*. Additionally, if a set $I \subseteq E$ is a subset of any base, then $I$ is called *independent*. The definition above generalizes the classical notion of matroids [9, 12], which can only have finite ground sets. Mathlib implements matroids as follows (this is a formalization of the definition above):

```
def Matroid.ExchangeProperty {α : Type*} (P : Set α → Prop) : Prop :=
  ∀ X Y : Set α, P X → P Y → ∀ a ∈ X \ Y, ∃ b ∈ Y \ X, P (insert b (X \ {a}))

def Maximal (P : α → Prop) (x : α) : Prop :=
  P x ∧ ∀ y : α, P y → x ≤ y → y ≤ x

def Matroid.ExistsMaximalSubsetProperty {α : Type*} (P : Set α → Prop) (X : Set α) : Prop :=
  ∀ I : Set α, P I → I ⊆ X →
    ∃ J : Set α, I ⊆ J ∧ Maximal (fun K : Set α => P K ∧ K ⊆ X) J

structure Matroid (α : Type*) where
  (E : Set α)
  (IsBase : Set α → Prop)
  (Indep : Set α → Prop)
  (indep_iff' : ∀ I : Set α, Indep I ↔ ∃ B : Set α, IsBase B ∧ I ⊆ B)
  (exists_isBase : ∃ B : Set α, IsBase B)
  (isBase_exchange : Matroid.ExchangeProperty IsBase)
  (maximality : ∀ X : Set α, X ⊆ E → Matroid.ExistsMaximalSubsetProperty Indep X)
  (subset_ground : ∀ B : Set α, IsBase B → B ⊆ E)
```

Additionally, Mathlib allows the user to construct matroids (potentially infinite) in terms of the *independence conditions* using:

```
structure IndepMatroid (α : Type*) where
  (E : Set α)
  (Indep : Set α → Prop)
  (indep_empty : Indep ∅)
  (indep_subset : ∀ I J : Set α, Indep J → I ⊆ J → Indep I)
  (indep_aug : ∀ I B : Set α, Indep I →
    ¬ Maximal Indep I → Maximal Indep B → ∃ x ∈ B \ I, Indep (insert x I))
  (indep_maximal : ∀ X : Set α, X ⊆ E → Matroid.ExistsMaximalSubsetProperty Indep X)
  (subset_ground : ∀ I : Set α, Indep I → I ⊆ E)
```

One can then obtain `Matroid α` from `IndepMatroid α` via `IndepMatroid.matroid`. The independence conditions frequently appear in constructions and proofs in classical literature [9, 12], and we use `IndepMatroid` to construct matroids in our library.

Though we generally work with infinite matroids, our final results require that the matroids have finite rank. A *finite-rank* matroid is one that has a finite base, defined in Mathlib as follows:

```
class RankFinite {α : Type*} (M : Matroid α) : Prop where
  exists_finite_isBase : ∃ B : Set α, M.IsBase B ∧ B.Finite
```

### 3.2 Totally Unimodular Matrices

In our work, regular matroids are defined in terms of totally unimodular matrices [9, 12]. Before introducing their definition, let us review how matrices and submatrices are implemented in Mathlib. A matrix with rows indexed by `m`, columns indexed by `n`, and entries of type $\alpha$ is represented by `Matrix m n α`, implemented as a (curried [11]) binary function `m → n → α`. Thus, the elements of matrix `A` can be accessed with `A i j`. Similarly, `Matrix.submatrix` is defined so that `(A.submatrix f g) i j = A (f i) (g j)` holds. Note that `Matrix.submatrix` may repeat and reorder rows and columns. For example, if

$$A = \begin{bmatrix} 1 & 2 & 3 \\ 4 & 5 & 6 \\ 7 & 8 & 9 \end{bmatrix}, \quad \texttt{f = ![0]}, \quad \texttt{g = ![2, 2, 0, 0]},$$

then `A.submatrix f g` $= \begin{bmatrix} 3 & 3 & 1 & 1 \end{bmatrix}$, typed as a matrix, not a vector.

Now, a matrix $A$ over a commutative ring $R$ is called *totally unimodular* if every finite square submatrix of $A$ (not necessarily contiguous, with no row or column taken twice) has determinant in $\{-1, 0, 1\}$. Mathlib implements this definition as follows:

```
def Matrix.IsTotallyUnimodular {m n R : Type*} [CommRing R] (A : Matrix m n R) : Prop :=
  ∀ k : ℕ, ∀ f : Fin k → m, ∀ g : Fin k → n, f.Injective → g.Injective →
      (A.submatrix f g).det ∈ Set.range SignType.cast
```

Here, `SignType` is an inductive type with three values: `zero`, `neg`, and `pos`; and `SignType.cast` maps them to `(0 : R)`, `(-1 : R)`, and `(1 : R)`, respectively.

Note that the indexing functions `f` and `g` are required to be injective in the definition, but this condition can be lifted:

```
lemma Matrix.isTotallyUnimodular_iff {m n R : Type*} [CommRing R] (A : Matrix m n R) :
    A.IsTotallyUnimodular ↔
    ∀ k : ℕ, ∀ f : Fin k → m, ∀ g : Fin k → n,
      (A.submatrix f g).det ∈ Set.range SignType.cast
```

Similarly, lemma `Matrix.isTotallyUnimodular_iff_fintype` equivalently characterizes total unimodularity by quantifying over any `Fintype` from any universe in place of `Fin k` above.

Keep in mind that the determinant is computed over $R$, so for certain commutative rings, all matrices are trivially totally unimodular, for example, for $R = \mathbb{Z}_3$.

### 3.3 Types and Subsets

In our project, we often have the following terms in the context:

```
(α : Type) (E : Set α) (I : Set α) (hIE : I ⊆ E)
```

Depending on the situation, there are three ways we may treat the set `I`. First, it may be viewed as a set of elements of type $\alpha$, its original type, so we simply write `I`. Second, we may need to re-type `I` as a set of elements of the type `E.Elem`. Then we write `E ↓∩ I` using notation from Mathlib. Finally, `I` may be used as a set of elements of the type `I.Elem`. In this case, we write `Set.univ` of the correct type, which is usually inferred from the context.

### 3.4 Block Matrices

In this project, we often construct matrices by composing them from blocks using the following Mathlib definitions:

- `Matrix.fromRows A₁ A₂` constructs $\begin{array}{|c|}\hline A_1 \\ \hline A_2 \\ \hline \end{array}$
- `Matrix.fromCols A₁ A₂` constructs $\begin{array}{|c|c|}\hline A_1 & A_2 \\ \hline \end{array}$
- `Matrix.fromBlocks A₁₁ A₁₂ A₂₁ A₂₂` constructs $\begin{array}{|c|c|}\hline A_{11} & A_{12} \\ \hline A_{21} & A_{22} \\ \hline \end{array}$

## 4 Re-typing Matrix Dimensions

When constructing matroids, we often need to convert a block matrix whose blocks are indexed by disjoint sets into a matrix indexed by unions of those index sets. Although the contents of the matrix stay the same, both its dimensions change their type from a `Sum` of sets to a `Set` union of those sets. To this end, we implemented:

```
def Subtype.toSum {α : Type*} {X Y : Set α}
    [∀ a, Decidable (a ∈ X)]
    [∀ a, Decidable (a ∈ Y)]
    (i : (X ∪ Y).Elem) :
    X.Elem ⊕ Y.Elem :=
  if hiX : i.val ∈ X then Sum.inl ⟨i, hiX⟩ else
  if hiY : i.val ∈ Y then Sum.inr ⟨i, hiY⟩ else
  (i.property.elim hiX hiY).elim
```

This allows us to re-type matrix dimensions and thus define the matrix transformation `Matrix.toMatrixUnionUnion` so that `A.toMatrixUnionUnion i j = A i.toSum j.toSum` holds. We also define an auxiliary function `Matrix.toMatrixElemElem` for convenience, but it is not a part of the trusted code.

## 5 Vector Matroids

Vector matroids [9, 12] is the most fundamental matroid class formalized in our work, serving as the basis for binary and regular matroids in later sections. A *vector matroid* is constructed from a matrix $A$ by taking the column index set as the ground set and declaring a set $I$ to be independent if the set of columns of $A$ indexed by $I$ is linearly independent. To capture this theoretical definition, we first implement the independence predicate as:

```
def Matrix.IndepCols {α R : Type*} {X Y : Set α} [Semiring R]
    (A : Matrix X Y R) (I : Set α) :
    Prop :=
  I ⊆ Y ∧ LinearIndepOn R Aᵀ (Y ↓∩ I)
```

Next, we construct an `IndepMatroid`:

```
def Matrix.toIndepMatroid {α R : Type*} {X Y : Set α} [DivisionRing R]
    (A : Matrix X Y R) :
    IndepMatroid α where
  E := Y
  Indep := A.IndepCols
  indep_empty := A.indepCols_empty
  indep_subset := A.indepCols_subset
  indep_aug := A.indepCols_aug
  indep_maximal S _ := A.indepCols_maximal S
  subset_ground _ := And.left
```

Finally, we convert `IndepMatroid` to `Matroid`:

```
def Matrix.toMatroid {α R : Type*} {X Y : Set α} [DivisionRing R]
    (A : Matrix X Y R) :
    Matroid α :=
  A.toIndepMatroid.matroid
```

Going forward, we use `Matrix.toMatroid` for constructing vector matroids from matrices.

As part of the construction above, we had to show that `Matrix.IndepCols` satisfies the so-called *augmentation property*: if $I$ is a non-maximal independent set and $J$ is a maximal independent set, then there exists an element $x \in J \setminus I$ such that $I \cup \{x\}$ is independent. It is worth noting that while we define `Matrix.IndepCols` over a semiring $R$ for the sake of generality, the augmentation property requires $R$ to be at least a division ring. Indeed, let $R = \mathbb{Z}_6$, which is in fact a ring, and consider

$$A = \begin{bmatrix} 0 & 1 & 2 & 3 \\ 1 & 0 & 3 & 2 \end{bmatrix}$$

with columns indexed by $\{0, 1, 2, 3\}$. Then $I = \{0\}$ is a non-maximal independent set and $J = \{2, 3\}$ is a maximal independent set over $R$, but they do not satisfy the augmentation property. For this reason, we require $R$ to be a division ring in the augmentation property and all dependent results.

Additionally, we show that vector matroids as defined above are finitary, i.e., an infinite subset in a vector matroid is independent if and only if so are all its finite subsets:

```
lemma Matrix.toMatroid_isFinitary {α R : Type*} {X Y : Set α} [DivisionRing R]
    (A : Matrix X Y R) :
    A.toMatroid.Finitary
```

## 6 Standard Representations

The *standard representation* [9, 12] of a vector matroid is the following structure:

```
structure StandardRepr (α R : Type*) [DecidableEq α] where
  X : Set α
  Y : Set α
  hXY : Disjoint X Y
  B : Matrix X Y R
  decmemX : ∀ a, Decidable (a ∈ X)
  decmemY : ∀ a, Decidable (a ∈ Y)
```

In essence, this is a wrapper for the standard representation matrix $B$ indexed by disjoint sets $X$ and $Y$, bundled together with the membership decidability for $X$ and $Y$. The standard representation matrix $B$ corresponds to the full representation matrix $\begin{bmatrix} \mathbb{1} & B \end{bmatrix}$ with the conversion implemented as:

```
def StandardRepr.toFull {α R : Type*} [DecidableEq α] [Zero R] [One R]
    (S : StandardRepr α R) :
    Matrix S.X.Elem (S.X ∪ S.Y).Elem R :=
  ((Matrix.fromCols 1 S.B) · ∘ Subtype.toSum)
```

Thus, the vector matroid given by its standard representation is constructed as follows:

```
def StandardRepr.toMatroid {α R : Type*} [DecidableEq α] [DivisionRing R]
    (S : StandardRepr α R) :
    Matroid α :=
  S.toFull.toMatroid
```

In this matroid, the ground set is $X \cup Y$, and a set $I \subseteq X \cup Y$ is independent if the columns of $\begin{bmatrix} \mathbb{1} & B \end{bmatrix}$ indexed by $I$ are linearly independent over $R$.

Below are several results we prove about standard representations, which are either used in the proof of regularity of 1-, 2-, and 3-sums, or could be useful for downstream projects.

First, we show that if the row index set $X$ of a standard representation is finite, then $X$ is a base in the resulting matroid:

```
lemma StandardRepr.toMatroid_isBase_X {α R : Type*} [DecidableEq α] [Field R]
    (S : StandardRepr α R) [Fintype S.X] :
    S.toMatroid.IsBase S.X
```

This lemma characterizes what sets can serve as row index sets of standard representations and motivates the corresponding hypotheses in the code snippets below.

Next, we prove that a full representation of a vector matroid can be transformed into a standard representation of the same matroid, with a given base as the row index set:

```
lemma Matrix.exists_standardRepr_isBase {α R : Type*} [DecidableEq α] [DivisionRing R]
    {X Y G : Set α} (A : Matrix X Y R) (hAG : A.toMatroid.IsBase G) :
    ∃ S : StandardRepr α R, S.X = G ∧ S.toMatroid = A.toMatroid
```

In classical literature on matroid theory [9, 12], this follows by simply performing a sequence of elementary row operations akin to Gaussian elimination. Our formal proof used a different approach, utilizing Mathlib's results about bases and linear independence. First, we showed that the columns indexed by `G` form a basis of the module generated by all columns of `A`. Then we proved that performing a basis exchange yields the correct standard representation matrix.

We also prove an analog of the above lemma that additionally preserves total unimodularity of the representation matrix:

```
lemma Matrix.exists_standardRepr_isBase_isTotallyUnimodular {α R : Type*}
    [DecidableEq α] [Field R] {X Y G : Set α} [Fintype G]
    (A : Matrix X Y R) (hAG : A.toMatroid.IsBase G) (hA : A.IsTotallyUnimodular) :
    ∃ S : StandardRepr α R, S.X = G ∧ S.toMatroid = A.toMatroid ∧ S.B.IsTotallyUnimodular
```

Classical literature [9, 12] observes that elementary row operations preserve total unimodularity and then simply refers to the proof of the previous lemma. Unfortunately, we could not take advantage of such a reduction, as it would be hard to verify that total unimodularity is preserved in our prior approach. Thus, we implemented an inductive proof essentially following the ideas of [9, 12]. Note that this lemma takes stronger assumptions than the previous one, namely that `G` has to be finite and that multiplication in `R` has to commute.

Another result we prove is that two standard representations of the same vector matroid over $\mathbb{Z}_2$ with the same finite row index set must be identical:

```
lemma ext_standardRepr_of_same_matroid_same_X {α : Type*} [DecidableEq α]
    {S₁ S₂ : StandardRepr α Z2} [Fintype S₁.X]
    (hSS : S₁.toMatroid = S₂.toMatroid) (hXX : S₁.X = S₂.X) :
    S₁ = S₂
```

Although this particular lemma is not employed later in our project, it captures an important result that a binary matroid has an essentially unique standard representation [9, 12]. Nevertheless, we make use of a very similar result:

```
lemma support_eq_support_of_same_matroid_same_X {F₁ : Type u₁} {F₂ : Type u₂}
    {α : Type max u₁ u₂ v} [DecidableEq α]
    [DecidableEq F₁] [DecidableEq F₂] [Field F₁] [Field F₂]
    {S₁ : StandardRepr α F₁} {S₂ : StandardRepr α F₂} [Fintype S₂.X]
    (hSS : S₁.toMatroid = S₂.toMatroid) (hXX : S₁.X = S₂.X) :
    let hYY : S₁.Y = S₂.Y := sorry
    hXX ▸ hYY ▸ S₁.B.support = S₂.B.support
```

This states that two standard representations of a vector matroid with identical (finite) row index sets have the same support, i.e., the zeros in them appear on identical positions. Crucially, this holds for any two standard representations over any two fields (where equality is decidable), and we later use it for $\mathbb{Q}$ and $\mathbb{Z}_2$.

## 7 Regular Matroids

Regular matroids [9, 12] are the core subject of Seymour's theorem. A matroid is *regular* if it can be constructed (as a vector matroid) from a rational totally unimodular matrix:

```
def Matroid.IsRegular {α : Type*} (M : Matroid α) : Prop :=
  ∃ X Y : Set α, ∃ A : Matrix X Y ℚ, A.IsTotallyUnimodular ∧ A.toMatroid = M
```

One key result we prove is that every regular matroid is in fact *binary*, i.e., can be constructed from a binary matrix:

```
lemma Matroid.IsRegular.isBinary {α : Type*} [DecidableEq α] {M : Matroid α}
    (hM : M.IsRegular) :
    ∃ X : Set α, ∃ Y : Set α, ∃ A : Matrix X Y Z2, A.toMatroid = M
```

Another important lemma we prove about regular matroids is their equivalent characterization in terms of totally unimodular signings. First, let us introduce the necessary definitions. We say that a matrix is a *signing* of another matrix if their values are identical up to signs:

```
def Matrix.IsSigningOf {X Y R : Type*} [LinearOrderedRing R] {n : ℕ}
    (A : Matrix X Y R) (U : Matrix X Y (ZMod n)) :
    Prop :=
  ∀ i : X, ∀ j : Y, |A i j| = (U i j).val
```

We then say that a binary matrix *has a totally unimodular signing* if it has a signing matrix that is rational and totally unimodular:

```
def Matrix.IsTuSigningOf {X Y : Type*} (A : Matrix X Y ℚ) (U : Matrix X Y Z2) : Prop :=
  A.IsTotallyUnimodular ∧ A.IsSigningOf U

def Matrix.HasTuSigning {X Y : Type*} (U : Matrix X Y Z2) : Prop :=
  ∃ A : Matrix X Y ℚ, A.IsTuSigningOf U
```

Now, we can state the characterization: given a standard representation over $\mathbb{Z}_2$, its matrix has a totally unimodular signing if and only if the matroid obtained from the representation is regular.

```
lemma StandardRepr.toMatroid_isRegular_iff_hasTuSigning {α : Type*} [DecidableEq α]
    (S : StandardRepr α Z2) [Finite S.X] :
  S.toMatroid.IsRegular ↔ S.B.HasTuSigning
```

Out of all definitions in this section, only `Matroid.IsRegular` is a part of the trusted code.

## 8 The 1-Sum

All matroid sums are defined on three levels: the `Matrix` level, the `StandardRepr` level, and the `Matroid` level. Let us review the distribution of responsibilities between the three levels.

```
def matrixSum1 {R : Type*} [Zero R] {Xₗ Yₗ Xᵣ Yᵣ : Type*}
    (Aₗ : Matrix Xₗ Yₗ R) (Aᵣ : Matrix Xᵣ Yᵣ R) :
    Matrix (Xₗ ⊕ Xᵣ) (Yₗ ⊕ Yᵣ) R :=
  Matrix.fromBlocks Aₗ 0 0 Aᵣ
```

The same matrix in a picture:

| | |
|---|---|
| $A_\ell$ | 0 |
| 0 | $A_r$ |

The `Matrix` level defines the standard representation matrix of the output matroid as a matrix indexed by `Sum` of indexing types. This definition is so straightforward that it would be natural to inline it into the subsequent definition. However, we retained it as a separate declaration for consistency with the 2- and 3-sums, whose matrix constructions are more elaborate.

```
noncomputable def standardReprSum1 {α : Type*} [DecidableEq α]
    {Sₗ Sᵣ : StandardRepr α Z2}
    (hXY : Disjoint Sₗ.X Sᵣ.Y)
    (hYX : Disjoint Sₗ.Y Sᵣ.X) :
    Option (StandardRepr α Z2) :=
  open scoped Classical in if
    Disjoint Sₗ.X Sᵣ.X ∧ Disjoint Sₗ.Y Sᵣ.Y
  then
    some ⟨
      Sₗ.X ∪ Sᵣ.X,
      Sₗ.Y ∪ Sᵣ.Y,
      sorry,
      (matrixSum1 Sₗ.B Sᵣ.B).toMatrixUnionUnion,
      inferInstance,
      inferInstance⟩
  else
    none
```

The `StandardRepr` level builds on top of the `Matrix` level. It converts the output matrix from being indexed by `Sum` to being index by set unions, it provides a proof that the resulting standard representation again has row indices and column indices disjoint, and it checks whether the operation is valid — if the preconditions are not met, it outputs `none` instead of `some` standard representation.

```
def Matroid.IsSum1of {α : Type*} [DecidableEq α] (M : Matroid α) (Mₗ Mᵣ : Matroid α) :
    Prop :=
  ∃ S Sₗ Sᵣ : StandardRepr α Z2,
  ∃ hXY : Disjoint Sₗ.X Sᵣ.Y,
  ∃ hYX : Disjoint Sₗ.Y Sᵣ.X,
  standardReprSum1 hXY hYX = some S
  ∧ S.toMatroid = M
  ∧ Sₗ.toMatroid = Mₗ
  ∧ Sᵣ.toMatroid = Mᵣ
```

The `Matroid` level builds on top of the standard representation level but talks about matroids, the combinatorial objects. On the `Matroid` level, we do not define a function; instead, we define a predicate — when $M$ is a 1-sum of $M_\ell$ and $M_r$.

In addition to basic API about the 1-sum, we also provide a theorem `Matroid.IsSum1of.eq_disjointSum` that establishes the equality between the disjoint sum (defined in Mathlib) and the 1-sum (defined in our project) of binary matroids.

## 9 The 2-Sum

The definition on 2-sum is also implemented on the three levels.

```
def matrixSum2 {R : Type*} [Semiring R] {Xₗ Yₗ Xᵣ Yᵣ : Type*}
    (Aₗ : Matrix Xₗ Yₗ R) (r : Yₗ → R) (Aᵣ : Matrix Xᵣ Yᵣ R) (c : Xᵣ → R) :
    Matrix (Xₗ ⊕ Xᵣ) (Yₗ ⊕ Yᵣ) R :=
  Matrix.fromBlocks Aₗ 0 (fun i j => c i * r j) Aᵣ
```

The Matrix level is pretty similar to the one of the 1-sum. Again, the two given matrices are placed along the main diagonal of the resulting block matrix. The resulting two blocks are not both zero, however, as this time the bottom left matrix contains the outer product of two given vectors. The same matrix in a picture:

| | |
|---|---|
| $A_\ell$ | $0$ |
| $c \otimes r$ | $A_r$ |

```
noncomputable def standardReprSum2 {α : Type*} [DecidableEq α]
    {Sₗ Sᵣ : StandardRepr α Z2} {x y : α}
    (hXX : Sₗ.X ∩ Sᵣ.X = {x})
    (hYY : Sₗ.Y ∩ Sᵣ.Y = {y})
    (hXY : Disjoint Sₗ.X Sᵣ.Y)
    (hYX : Disjoint Sₗ.Y Sᵣ.X) :
    Option (StandardRepr α Z2) :=
  let Aₗ : Matrix (Sₗ.X \ {x}).Elem Sₗ.Y Z2 := Sₗ.B.submatrix Set.diff_subset.elem id
  let Aᵣ : Matrix Sᵣ.X (Sᵣ.Y \ {y}).Elem Z2 := Sᵣ.B.submatrix id Set.diff_subset.elem
  let r : Sₗ.Y.Elem → Z2 := Sₗ.B ⟨x, sorry⟩
  let c : Sᵣ.X.Elem → Z2 := (Sᵣ.B · ⟨y, sorry⟩)
  open scoped Classical in if
    r ≠ 0 ∧ c ≠ 0
  then
    some ⟨
      (Sₗ.X \ {x}) ∪ Sᵣ.X,
      Sₗ.Y ∪ (Sᵣ.Y \ {y}),
      sorry,
      (matrixSum2 Aₗ r Aᵣ c).toMatrixUnionUnion,
      inferInstance,
      inferInstance⟩
  else
    none
```

The StandardRepr level is more complicated. We first need to slice the last row of the matrix $S_\ell$.B and the first column of the matrix $S_r$.B as the two separate vectors (r and c), naming the two remaining matrices $A_\ell$ and $A_r$ respectively. To identify the special row and the special column (remember that matrices do not have rows and columns ordered, as much as we like to draw certain canonical ordering or rows and columns on paper for helpful visuals), we need to be given a specific element x in $S_\ell$.X ∩ $S_r$.X and a specific element y in $S_\ell$.Y ∩ $S_r$.Y and promised that there is no other element in any pairwise intersection among the four indexing sets. The following picture shows how $S_\ell$.B and $S_r$.B are taken apart:

$$S_\ell.\mathrm{B} = \begin{array}{|c|}\hline A_\ell \\ \hline r \\ \hline\end{array}\ , \quad S_r.\mathrm{B} = \begin{array}{|c|c|}\hline c & A_r \\ \hline\end{array}$$

Now we know the arguments to be given to the Matrix level. Again, we convert the output matrix from being indexed by Sum to being index by set unions, we provide a proof that the resulting standard representation has row indices and column indices disjoint, and we check whether the operation is valid — this time, the condition is that neither r nor c is a zero vector.

```
def Matroid.IsSum2of {α : Type*} [DecidableEq α] (M : Matroid α) (Mₗ Mᵣ : Matroid α) :
    Prop :=
  ∃ S Sₗ Sᵣ : StandardRepr α Z2,
  ∃ x y : α,
  ∃ hXX : Sₗ.X ∩ Sᵣ.X = {x},
  ∃ hYY : Sₗ.Y ∩ Sᵣ.Y = {y},
  ∃ hXY : Disjoint Sₗ.X Sᵣ.Y,
  ∃ hYX : Disjoint Sₗ.Y Sᵣ.X,
  standardReprSum2 hXX hYY hXY hYX = some S
  ∧ S.toMatroid = M
  ∧ Sₗ.toMatroid = Mₗ
  ∧ Sᵣ.toMatroid = Mᵣ
```

The `Matroid` level is again a predicate — when $M$ is a 2-sum of $M_\ell$ and $M_r$.

## 10 The 3-Sum

The 3-sum of binary matroids is defined as follows. Let $X_\ell$, $Y_\ell$, $X_r$, and $Y_r$ be sets with the following properties:

- $X_\ell \cap X_r = \{x_2, x_1, x_0\}$ for some distinct $x_0$, $x_1$, and $x_2$
- $Y_\ell \cap Y_r = \{y_0, y_1, y_2\}$ for some distinct $y_0$, $y_1$, and $y_2$
- $X_\ell \cap Y_\ell = X_\ell \cap Y_r = X_r \cap Y_\ell = X_r \cap Y_r = \emptyset$

Let $B_\ell \in \mathbb{Z}_2^{X_\ell \times Y_\ell}$ and $B_r \in \mathbb{Z}_2^{X_r \times Y_r}$ be matrices of the form

$$B_\ell = \left[\begin{array}{cc|c} \multicolumn{2}{c|}{A_\ell} & 0 \\ & \begin{array}{|c|c|}\hline 1 & 1 \\ \hline\end{array} & 0 \\ \hline D_\ell & D_0 & 1 \\ & & 1 \end{array}\right], \quad B_r = \left[\begin{array}{c|c|c} \begin{array}{c|c}1 & 1\end{array} & 0 & 0 \\ \hline D_0 & \begin{array}{c}1 \\ \hline 1\end{array} & A_r \\ \cline{1-2} D_r & & \end{array}\right]$$

where $D_0$ is invertible. Then their 3-sum is

$$B = \left[\begin{array}{cc|c|c} \multicolumn{2}{c|}{A_\ell} & \multicolumn{2}{c}{0} \\ & \begin{array}{c|c}1 & 1\end{array} & 0 & \\ \hline D_\ell & D_0 & \begin{array}{c}1 \\ \hline 1\end{array} & A_r \\ \cline{1-3} D_{\ell r} & D_r & & \end{array}\right] \quad \text{where } D_{\ell r} = D_r \cdot D_0^{-1} \cdot D_\ell$$

Here $D_0 \in \mathbb{Z}_2^{\{x_0,x_1\}\times\{y_0,y_1\}}$, $\left[\begin{array}{c|c} \begin{array}{c|c}1 & 1\end{array} & 0 \\ \hline D_0 & \begin{array}{c}1 \\ \hline 1\end{array} \end{array}\right] \in \mathbb{Z}_2^{\{x_2,x_0,x_1\}\times\{y_0,y_1,y_2\}}$, and the indexing is kept consistent between $B_\ell$, $B_r$, and $B$. Consequently, a matroid $M$ is a 3-sum of matroids $M_\ell$ and $M_r$ if they admit standard representations over $\mathbb{Z}_2$ with matrices $B$, $B_\ell$, and $B_r$ of the form above.

In our implementation, we frequently deal with sets with one, two, or three elements removed. To make our code more compact, we added abbreviations for removing one, two, and three elements from a set, as well as a definition for re-typing an element of a set with three elements removed as an element of the original set:

```
abbrev Set.drop1 {α : Type*} (Z : Set α) (z₀ : Z) : Set α :=
  Z \ {z₀.val}

abbrev Set.drop2 {α : Type*} (Z : Set α) (z₀ z₁ : Z) : Set α :=
  Z \ {z₀.val, z₁.val}

abbrev Set.drop3 {α : Type*} (Z : Set α) (z₀ z₁ z₂ : Z) : Set α :=
  Z \ {z₀.val, z₁.val, z₂.val}

def undrop3 {α : Type*} {Z : Set α} {z₀ z₁ z₂ : Z} (i : Z.drop3 z₀ z₁ z₂) : Z :=
  ⟨i.val, i.property.left⟩
```

Now, to define the 3-sum of matrices, we introduce a structure comprising the blocks of the summands:

```
structure MatrixSum3 (Xₗ Yₗ Xᵣ Yᵣ R : Type*) where
  Aₗ : Matrix (Xₗ ⊕ Unit) (Yₗ ⊕ Fin 2) R
  Dₗ : Matrix (Fin 2) Yₗ R
  D₀ₗ : Matrix (Fin 2) (Fin 2) R
  D₀ᵣ : Matrix (Fin 2) (Fin 2) R
  Dᵣ : Matrix Xᵣ (Fin 2) R
  Aᵣ : Matrix (Fin 2 ⊕ Xᵣ) (Unit ⊕ Yᵣ) R
```

Here D₀ₗ and D₀ᵣ refer to block $D_0$ in $B_\ell$ and $B_r$, respectively. It is then straightforward to define the resulting 3-sum matrix:

```
noncomputable def MatrixSum3.matrix {Xₗ Yₗ Xᵣ Yᵣ R : Type*} [CommRing R]
    (S : MatrixSum3 Xₗ Yₗ Xᵣ Yᵣ R) :
    Matrix ((Xₗ ⊕ Unit) ⊕ (Fin 2 ⊕ Xᵣ)) ((Yₗ ⊕ Fin 2) ⊕ (Unit ⊕ Yᵣ)) R :=
  Matrix.fromBlocks S.Aₗ 0
    (Matrix.fromBlocks S.Dₗ S.D₀ₗ
      (S.Dᵣ * S.D₀ₗ⁻¹ * S.Dₗ) S.Dᵣ
    ) S.Aᵣ
```

Introducing these definitions creates an abstraction layer that allows us to work with the blocks used to construct a 3-sum of matrices without the need to manually obtain them from the summands each time. Moreover, this drastically simplifies the implementation of results that require additional assumptions on the summands. Without these definitions, one has to repeatedly extract the blocks from the summands before the additional assumptions or the final result can be stated and in the proof as well, which is extremely cumbersome.

To further facilitate our implementation of the 3-sum, we pack the inner workings of obtaining the blocks from the summands into the following definition:

```
def blocksToMatrixSum3 {Xₗ Yₗ Xᵣ Yᵣ R : Type*}
    (Bₗ : Matrix ((Xₗ ⊕ Unit) ⊕ Fin 2) ((Yₗ ⊕ Fin 2) ⊕ Unit) R)
    (Bᵣ : Matrix (Unit ⊕ (Fin 2 ⊕ Xᵣ)) (Fin 2 ⊕ (Unit ⊕ Yᵣ)) R) :
    MatrixSum3 Xₗ Yₗ Xᵣ Yᵣ R where
  Aₗ := Bₗ.toBlocks₁₁
  Dₗ := Bₗ.toBlocks₂₁.toCols₁
  D₀ₗ := Bₗ.toBlocks₂₁.toCols₂
  D₀ᵣ := Bᵣ.toBlocks₂₁.toRows₁
  Dᵣ := Bᵣ.toBlocks₂₁.toRows₂
  Aᵣ := Bᵣ.toBlocks₂₂
```

This definition is particularly compact thanks to us changing the types of dimensions of $B_\ell$ and $B_r$. The corresponding transformation of dimensions of $B_\ell$ is then implemented as:

```
def Matrix.toBlockSummandₗ {α R : Type*} {Xₗ Yₗ : Set α} (Bₗ : Matrix Xₗ Yₗ R)
    (x₀ x₁ x₂ : Xₗ) (y₀ y₁ y₂ : Yₗ) :
    Matrix ((Xₗ.drop3 x₀ x₁ x₂ ⊕ Unit) ⊕ Fin 2) ((Yₗ.drop3 y₀ y₁ y₂ ⊕ Fin 2) ⊕ Unit) R :=
  Bₗ.submatrix
    (·.casesOn (·.casesOn undrop3 ↓x₂) ![x₀, x₁])
    (·.casesOn (·.casesOn undrop3 ![y₀, y₁]) ↓y₂)
```

We re-index $B_r$ analogously via `Matrix.toBlockSummandᵣ`.

Now, to implement 3-sums of standard representations, we perform one last reindexing to transform the dimensions of `MatrixSum3.matrix` into unions of sets via:

```
def Matrix.toMatrixDropUnionDrop {α : Type*} [DecidableEq α] {Xₗ Yₗ Xᵣ Yᵣ : Set α} {R : Type*}
    [∀ a, Decidable (a ∈ Xₗ)]
    [∀ a, Decidable (a ∈ Yₗ)]
    [∀ a, Decidable (a ∈ Xᵣ)]
    [∀ a, Decidable (a ∈ Yᵣ)]
    {x₀ₗ x₁ₗ x₂ₗ : Xₗ} {y₀ₗ y₁ₗ y₂ₗ : Yₗ}
    {x₀ᵣ x₁ᵣ x₂ᵣ : Xᵣ} {y₀ᵣ y₁ᵣ y₂ᵣ : Yᵣ}
    (A : Matrix
      ((Xₗ.drop3 x₀ₗ x₁ₗ x₂ₗ ⊕ Unit) ⊕ (Fin 2 ⊕ Xᵣ.drop3 x₀ᵣ x₁ᵣ x₂ᵣ))
      ((Yₗ.drop3 y₀ₗ y₁ₗ y₂ₗ ⊕ Fin 2) ⊕ (Unit ⊕ Yᵣ.drop3 y₀ᵣ y₁ᵣ y₂ᵣ))
      R) :
    Matrix (Xₗ.drop2 x₀ₗ x₁ₗ ∪ Xᵣ.drop1 x₂ᵣ).Elem (Yₗ.drop1 y₂ₗ ∪ Yᵣ.drop2 y₀ᵣ y₁ᵣ).Elem R :=
  A.submatrix
    (fun i : (Xₗ.drop2 x₀ₗ x₁ₗ ∪ Xᵣ.drop1 x₂ᵣ).Elem =>
      if hi₂ₗ : i.val = x₂ₗ then Sum.inl (Sum.inr 0) else
      if hiXₗ : i.val ∈ Xₗ.drop3 x₀ₗ x₁ₗ x₂ₗ then Sum.inl (Sum.inl ⟨i, hiXₗ⟩) else
      if hi₀ᵣ : i.val = x₀ᵣ then Sum.inr (Sum.inl 0) else
      if hi₁ᵣ : i.val = x₁ᵣ then Sum.inr (Sum.inl 1) else
      if hiXᵣ : i.val ∈ Xᵣ.drop3 x₀ᵣ x₁ᵣ x₂ᵣ then Sum.inr (Sum.inr ⟨i, hiXᵣ⟩) else
      False.elim sorry)
    (fun j : (Yₗ.drop1 y₂ₗ ∪ Yᵣ.drop2 y₀ᵣ y₁ᵣ).Elem =>
      if hj₀ₗ : j.val = y₀ₗ then Sum.inl (Sum.inr 0) else
      if hj₁ₗ : j.val = y₁ₗ then Sum.inl (Sum.inr 1) else
      if hjYₗ : j.val ∈ Yₗ.drop3 y₀ₗ y₁ₗ y₂ₗ then Sum.inl (Sum.inl ⟨j, hjYₗ⟩) else
      if hj₂ᵣ : j.val = y₂ᵣ then Sum.inr (Sum.inl 0) else
      if hjYᵣ : j.val ∈ Yᵣ.drop3 y₀ᵣ y₁ᵣ y₂ᵣ then Sum.inr (Sum.inr ⟨j, hjYᵣ⟩) else
      False.elim sorry)
```

This allows us to define the 3-sum of standard representations as follows:

```
noncomputable def standardReprSum3 {α : Type*} [DecidableEq α]
    {Sₗ Sᵣ : StandardRepr α Z2} {x₀ x₁ x₂ y₀ y₁ y₂ : α}
    (hXX : Sₗ.X ∩ Sᵣ.X = {x₀, x₁, x₂})
    (hYY : Sₗ.Y ∩ Sᵣ.Y = {y₀, y₁, y₂})
    (hXY : Disjoint Sₗ.X Sᵣ.Y)
    (hYX : Disjoint Sₗ.Y Sᵣ.X) :
    Option (StandardRepr α Z2) :=
  let x₀ₗ : Sₗ.X := ⟨x₀, sorry⟩
  let x₁ₗ : Sₗ.X := ⟨x₁, sorry⟩
  let x₂ₗ : Sₗ.X := ⟨x₂, sorry⟩
  let y₀ₗ : Sₗ.Y := ⟨y₀, sorry⟩
  let y₁ₗ : Sₗ.Y := ⟨y₁, sorry⟩
  let y₂ₗ : Sₗ.Y := ⟨y₂, sorry⟩
  let x₀ᵣ : Sᵣ.X := ⟨x₀, sorry⟩
  let x₁ᵣ : Sᵣ.X := ⟨x₁, sorry⟩
  let x₂ᵣ : Sᵣ.X := ⟨x₂, sorry⟩
  let y₀ᵣ : Sᵣ.Y := ⟨y₀, sorry⟩
  let y₁ᵣ : Sᵣ.Y := ⟨y₁, sorry⟩
  let y₂ᵣ : Sᵣ.Y := ⟨y₂, sorry⟩
  open scoped Classical in if
    ((x₀ ≠ x₁ ∧ x₀ ≠ x₂ ∧ x₁ ≠ x₂) ∧ (y₀ ≠ y₁ ∧ y₀ ≠ y₂ ∧ y₁ ≠ y₂))
    ∧ Sₗ.B.submatrix ![x₀ₗ, x₁ₗ] ![y₀ₗ, y₁ₗ] = Sᵣ.B.submatrix ![x₀ᵣ, x₁ᵣ] ![y₀ᵣ, y₁ᵣ]
    ∧ IsUnit (Sₗ.B.submatrix ![x₀ₗ, x₁ₗ] ![y₀ₗ, y₁ₗ])
    ∧ Sₗ.B x₀ₗ y₂ₗ = 1
    ∧ Sₗ.B x₁ₗ y₂ₗ = 1
```

```
    ∧ S₁.B x₂₁ y₀₁ = 1
    ∧ S₁.B x₂₁ y₁₁ = 1
    ∧ (∀ x : α, ∀ hx : x ∈ S₁.X, x ≠ x₀ ∧ x ≠ x₁ → S₁.B ⟨x, hx⟩ y₂₁ = 0)
    ∧ Sᵣ.B x₀ᵣ y₂ᵣ = 1
    ∧ Sᵣ.B x₁ᵣ y₂ᵣ = 1
    ∧ Sᵣ.B x₂ᵣ y₀ᵣ = 1
    ∧ Sᵣ.B x₂ᵣ y₁ᵣ = 1
    ∧ (∀ y : α, ∀ hy : y ∈ Sᵣ.Y, y ≠ y₀ ∧ y ≠ y₁ → Sᵣ.B x₂ᵣ ⟨y, hy⟩ = 0)
  then
    some ⟨
      (S₁.X.drop2 x₀₁ x₁₁) ∪ (Sᵣ.X.drop1 x₂ᵣ),
      (S₁.Y.drop1 y₂₁) ∪ (Sᵣ.Y.drop2 y₀ᵣ y₁ᵣ),
      sorry,
      (blocksToMatrixSum3
        (S₁.B.toBlockSummand₁ x₀₁ x₁₁ x₂₁ y₀₁ y₁₁ y₂₁)
        (Sᵣ.B.toBlockSummandᵣ x₀ᵣ x₁ᵣ x₂ᵣ y₀ᵣ y₁ᵣ y₂ᵣ)
      ).matrix.toMatrixDropUnionDrop,
      inferInstance,
      inferInstance⟩
  else
    none
```

The resulting definition has similar advantages to its analogs for the 1- and 2-sum:

- The data required to construct the 3-sum together with all intermediate objects and assumptions appear as named arguments.
- Conditions that are not needed to carry out the construction but necessary for the result to be valid are anonymous and appear in the `if` statement.
- The result is given by an `Option`, which evaluates to `some` standard representation if the produced 3-sum is valid, or `none` otherwise.

Finally, the `Matroid`-level predicate `Matroid.IsSum3of` is defined similarly to those for 1- and 2-sums by, ensuring consistency:

```
recall Matroid.IsSum3of {α : Type*} [DecidableEq α] (M M₁ Mᵣ : Matroid α) :
    Prop :=
  ∃ S S₁ Sᵣ : StandardRepr α Z2,
  ∃ x₀ x₁ x₂ y₀ y₁ y₂ : α,
  ∃ hXX : S₁.X ∩ Sᵣ.X = {x₀, x₁, x₂},
  ∃ hYY : S₁.Y ∩ Sᵣ.Y = {y₀, y₁, y₂},
  ∃ hXY : Disjoint S₁.X Sᵣ.Y,
  ∃ hYX : Disjoint S₁.Y Sᵣ.X,
  standardReprSum3 hXX hYY hXY hYX = some S
  ∧ S.toMatroid = M
  ∧ S₁.toMatroid = M₁
  ∧ Sᵣ.toMatroid = Mᵣ
```

## 11 Sums Preserve Regularity

In our library, the final theorems that regularity is preserved under 1-, 2-, and 3-sums are stated as follows.

```
theorem Matroid.IsSum1of.isRegular {α : Type*} [DecidableEq α] {M M₁ Mᵣ : Matroid α} :
  M.IsSum1of M₁ Mᵣ → M.RankFinite → M₁.IsRegular → Mᵣ.IsRegular → M.IsRegular

theorem Matroid.IsSum2of.isRegular {α : Type*} [DecidableEq α] {M M₁ Mᵣ : Matroid α} :
  M.IsSum2of M₁ Mᵣ → M.RankFinite → M₁.IsRegular → Mᵣ.IsRegular → M.IsRegular

theorem Matroid.IsSum3of.isRegular {α : Type*} [DecidableEq α] {M M₁ Mᵣ : Matroid α} :
  M.IsSum3of M₁ Mᵣ → M.RankFinite → M₁.IsRegular → Mᵣ.IsRegular → M.IsRegular
```

Note that these three theorems are stated for matroids and have the same interface. Moreover, when applying one of these results, a user is able to provide different representations for witnessing that $M$ is a 1-, 2-, or 3-sum of $M_\ell$ and $M_r$, for witnessing that $M$ has finite rank, and for witnessing that $M_\ell$ and $M_r$ are regular.

We split the proof of each of these theorems into three stages corresponding to the three abstraction layers used for the definitions: `Matroid`, `StandardRepr`, and `Matrix`.

The final `Matroid`-level theorems for all 1-, 2-, and 3-sums are reduced to the respective lemmas for standard representations by applying out lemmas `StandardRepr.toMatroid_isRegular_iff_hasTuSigning` and `StandardRepr.finite_X_of_toMatroid_rankFinite` in all three proofs. The reductions from the `StandardRepr` level to the `Matrix` level for 1- and 2-sums is straightforward — plug the standard representation matrices and their (rational) signings into `matrixSum1` and `matrixSum2`, respectively. For 3-sums, this reduction is more involved, as we additionally apply the following lemma to simplify the assumption on $D_0$:

```
lemma Matrix.isUnit_2x2 (A : Matrix (Fin 2) (Fin 2) Z2) (hA : IsUnit A) :
  ∃ f : Fin 2 ≃ Fin 2, ∃ g : Fin 2 ≃ Fin 2,
    A.submatrix f g = 1 ∨ A.submatrix f g = !![1, 1; 0, 1]
```

Therefore, up to reindexing, $D_0$ is either $\begin{bmatrix} 1 & 0 \\ 0 & 1 \end{bmatrix}$ or $\begin{bmatrix} 1 & 1 \\ 0 & 1 \end{bmatrix}$. Performing the reduction at this stage allows us to invoke `Matrix.isUnit_2x2` only once and then simply consider the two special forms of $D_0$.

On the `Matrix` level, our formal proof that 1-sums preserve total unimodularity of matrices is nearly identical to [12]. For 2-sums, we streamlined the proof by reformulating it as a forward argument by induction. For 3-sums, the entire argument was significantly reworked to simplify and streamline the approach of [12]. On a high level, we make two major changes, which we discuss in detail below.

The first key difference is that we re-sign the summands only once, rather than multiple times. Like in [12], we start with totally unimodular signings exhibiting regularity of the two summands. Then we multiply their rows and columns by $\pm 1$ factors (which preserves total unimodularity) so that the submatrix $\begin{array}{|c|c|c|}\hline 1 & 1 & 0 \\ \hline \multicolumn{2}{|c|}{\multirow{2}{*}{$D_0$}} & 1 \\ \cline{3-3} \multicolumn{2}{|c|}{} & 1 \\ \hline \end{array}$ is signed in both summands simultaneously as either $\begin{bmatrix} 1 & 1 & 0 \\ 1 & 0 & 1 \\ 0 & -1 & 1 \end{bmatrix}$ or $\begin{bmatrix} 1 & 1 & 0 \\ 1 & 1 & 1 \\ 0 & 1 & 1 \end{bmatrix}$, depending on whether $D_0$ is $\begin{bmatrix} 1 & 0 \\ 0 & 1 \end{bmatrix}$ or $\begin{bmatrix} 1 & 1 \\ 0 & 1 \end{bmatrix}$. Thus, we get totally unimodular signings of the summands that coincide on the intersection, which allows us to define the *canonical* signing of the entire 3-sum: use the same signs as in the re-signed summands everywhere except for the bottom-left block, which is signed via $D'_{\ell r} = D'_r \cdot (D'_0)^{-1} \cdot D'_\ell$, and the 0 block, which remains as is.

The main advantage of our approach is that we avoid chained constructions and proofs of properties of such constructions, and we do not need to define $\Delta$-sums. Moreover, unlike [12], our proof does not rely on the general lemma about re-signing totally unimodular matrices. This detail is crucial, as the proof of this lemma in [12] involves a graph-theoretic argument, which would be very challenging to formalize in Lean with the current tools available in Mathlib.

The second major difference from the approach of [12] is that our main argument does not deal with signings of 3-sums directly. Instead, we work with a matrix family called `MatrixLikeSum3` in our code. This allows us to split the proof of regularity of 3-sums into three clear steps. First, we show that pivoting on a non-zero entry in the top-left block of any matrix from this family produces a matrix that also belongs to this family. Next, we utilize the result from the first step to prove that every matrix in this family is totally unimodular. We do this via a similar argument to the proof that 2-sums of totally unimodular matrices are totally unimodular. Finally, we show that every canonical signing of the 3-sum matrix defined above is included in this matrix family and is thus totally unimodular. Overall, this proof takes a more systematic approach to deriving properties of signings of 3-sums and using them to prove their total unimodularity. Additionally, it conveniently reuses a large portion of the argument for 2-sums.

In some proofs, we worked with large case splits with up to 896 cases. To handle such situations, we used

```
all_goals try
```

followed by one or more tactics, discharging multiple goals at once without selecting them by hand or repeating the proof. We repeatedly applied this method to discharge the remaining goals in waves until the proof was complete.

## 12 Related Work

In Lean 4, the largest library formalizing matroid theory is due to Peter Nelson[3]. It implements infinite matroids following [3] together with many key notions and results about them. The definition that is fully formalized and is the most related to our work is `Matroid.disjointSum`. For binary matroids, this definition is equivalent to the 1-sum implemented in this paper. Moreover, it can be used for any matroids with disjoint ground sets, while our implementation is restricted to vector matroids constructed from $\mathbb{Z}_2$ matrices. Peter Nelson's repository also makes progress towards formalizing other related notions, such as representable matroids, though this work is still ongoing. It is also worth noting that the results in Mathlib[4] have been copied over from this repository and comprise a strict subset of it.

Building upon Peter Nelson's work, Gusakov's thesis [5] formalizes the proof of Tutte's excluded minor theorem and to this end implements definitions and results about representable matroids. The thesis formalizes representations and standard representations of matroids, which we also do in our work, but it takes a different approach. In particular, instead of working with matrix representations, the thesis implements a representation of `Matroid` $\alpha$ as a mapping from the entire type $\alpha$ to a vector space, which maps non-elements of the matroid to the zero vector and independent sets to linearly independent vectors. The advantage of this approach is that certain proofs become easier to formalize, but this comes at a cost of making it harder to match the implementation with the theory and believe the correctness of the code.

There are also two Lean 3 repositories due to Artem Vasilyev[5] and Bryan Gin-ge Chen[6] dedicated to formalization of matroid theory. Both of them work with finite matroids following [9] and implement basic definitions and properties of matroids concerning circuits, bases, and rank functions. These results are completely subsumed by the current implementation of matroids in Mathlib.

Jonas Keinholz [6] formalizes the classical definition of (finite) matroids [9, 12] in Isabelle/HOL along with other basic ideas such as minors, bases, circuits, rank, and closure. More recently, Wan et al. [13] use Keinholz's formalization to design a verification framework using a Locale that checks if a given collection of subsets of a given set is a matroid. The authors then showcase the verification algorithm by checking that the 0-1 knapsack problem does not conform to the matroid structure, while the fractional knapsack problem does. In comparison, Lean 4's Mathlib implements a more general definition of matroids and formalizes more results about them than either [6] or [13], but Lean lacks a procedure for formally verifying if a collection of sets has matroid structure.

In the HOL Light GitHub repository[7], John Harrison formalizes finitary matroids. The formalization closely follows the field theory notes of Pete L. Clark[8]. In particular, finitary matroids are defined in terms of a closure operator with similar properties as those proposed in [3]. This repository also includes a formal proof that this notion of (finitary) matroids is equivalent to the definition of a matroid using independent sets. Unlike Lean 4's Mathlib formalization (which includes formalizations of the closure operator and the notions of spanning sets), however, this notion of infinite matroids does not respect the notion of duality that is defined for matroids in [9, 12] as noted by [3].

Grzegorz Bancerek and Yasunari Shidama [1] formalize matroids in Mizar. Their formalization includes basic notions like rank, basis, and cycle as well as examples like the matroid of linearly independent subsets for a given vector space. Overall, the scope of the Mizar formalization is comparable to the Isabelle/HOL formalization, except that the Mizar formalization allows for infinite matroids. In this sense, it is comparable to the Lean definition in Mathlib, which also allows for infinite matroids. However, whereas Mizar uses independence conditions to define matroids, Lean uses base conditions for the main definition and provides an API for constructing matroids via independence conditions.

[3] `https://github.com/apnelson1/lean-matroids`

[4] `https://github.com/leanprover-community/mathlib4/tree/master/Mathlib/Combinatorics/Matroid`

[5] `https://github.com/VArtem/lean-matroids`

[6] `https://github.com/bryangingechen/lean-matroids`

[7] `https://github.com/jrh13/hol-light/blob/master/Library/matroids.ml`

[8] `https://plclark.github.io/PeteLClark/Expositions/FieldTheory.pdf`

## 13 Conclusion

In this work, we formally stated Seymour's decomposition theorem for regular matroids and implemented a formally verified proof of the forward (composition) direction of this theorem in the setting where the matroids have finite rank and may have infinite ground sets. To this end, we developed a modular and extensible library in Lean 4 formalizing definitions and lemmas about totally unimodular matrices, vector matroids, regular matroids, and 1-, 2-, and 3-sums of matrices, standard representations of vector matroids, and matroids. Our work demonstrates that one can effectively use Lean and Mathlib to formally verify advanced results from matroid theory and extend classical results to a more general setting.

The most natural continuation of our project is proving the decomposition direction of Seymour's theorem, stated as `Matroid.IsRegular.isGood` in our library. Our work can also serve as a starting point for formalizing Seymour's theorem for matroids of infinite rank [2].

## 14 Acknowledgments

We would like to thank Jasmin Blanchette for frequent consultations, Pietro Monticone for help with LeanProject, Damiano Testa for help with linters, Riccardo Brasca for a proof of `Matrix.one_linearIndependent`, Johan Commelin and Edward van de Meent for advice about proving `Matrix.fromBlocks_isTotallyUnimodular`, Aaron Liu for advice about handling `HEq`, Yaël Dillies for advice about inverting functions, and Jireh Loreaux for advice about singular matrices.

B.-H. Hwang was supported by a KIAS Individual Grant (MG098201) at Korea Institute for Advanced Study.

# A Suggestions for Lean's Ecosystem

When working on the project, our overall experience of working with Lean, Mathlib, and Aesop was very positive. Nevertheless, in our opinion, certain things could be made even better, and we share our improvement suggestions below.

## A.1 Suggestions for Lean

The `recall` command is really helpful for presentation and believability of our trusted code. Unfortunately, it currently does not support structures, which is a significant limitation. For example, suppose we want to recall

```
structure StandardRepr (α R : Type*) [DecidableEq α] where
  X : Set α
  Y : Set α
  hXY : Disjoint X Y
  B : Matrix X Y R
  decmemX : ∀ a, Decidable (a ∈ X)
  decmemY : ∀ a, Decidable (a ∈ Y)
```

Currently, we have to separately `recall` every field of the structure:

```
recall StandardRepr.X {α R : Type*} [DecidableEq α] :
  StandardRepr α R → Set α
recall StandardRepr.Y {α R : Type*} [DecidableEq α] :
  StandardRepr α R → Set α
recall StandardRepr.B {α R : Type*} [DecidableEq α] (S : StandardRepr α R) :
  Matrix S.X S.Y R
recall StandardRepr.hXY {α R : Type*} [DecidableEq α] (S : StandardRepr α R) :
  Disjoint S.X S.Y
recall StandardRepr.decmemX {α R : Type*} [DecidableEq α] (S : StandardRepr α R) :
  ∀ a, Decidable (a ∈ S.X)
recall StandardRepr.decmemY {α R : Type*} [DecidableEq α] (S : StandardRepr α R) :
  ∀ a, Decidable (a ∈ S.Y)
```

This code contains a lot of repetition and can be difficult to read, compared to the original declaration. Moreover, it does not ensure that all fields of the structure are listed, so we additionally check the constructor using `#guard_msgs`:

```
/--
info: StandardRep.mk.{u_1, u_2}
  {α : Type u_1} {R : Type u_2} [DecidableEq α]
  (X Y : Set α) (hXY : X ⫗ Y)
  (B : Matrix (X) (Y) R)
  (decmemX : (a : α) → Decidable (a ∈ X))
  (decmemY : (a : α) → Decidable (a ∈ Y)) :
  StandardRepr α R
-/
#guard_msgs in
#check StandardRepr.mk
```

Without this additional check, it would be possible to cheat. One could maliciously add a field claiming that $1 + 1 = 3$ to `StandardRepr`, not recall it in the presentation file, and then use it to prove any result without having to capture its mathematical essence. Although the approach of recalling all fields of a structure and then checking the constructor is sufficient for our presentation file, being able to `recall` the structure directly would make it easier to write complete and believable trusted code.

When refactoring our library, we found it tedious to identify definitions and lemmas that are unused and can be safely pruned. In many cases, we could identify them using project-wide search or the "go to reference" command in VS Code, but we had to do it by hand in every instance, and we had to be particularly careful whenever dot notation was used. Moreover, for lemmas with the `@[simp]` attribute, neither of the two approaches could be used to find the `simp` calls using such lemmas. Ultimately, the only way to know for sure if something was used or not was to try removing it and recompiling the entire project, which could take a long time especially if the change affected a file close to the root of the dependency graph. It would be ideal to have a tool that could show which definitions and lemmas are used in the implementation of each statement and proof, perhaps in the form of a dependency graph similar to the one generated by `leanblueprint`. This would not only help identify code that is never used, but also give a clear and precise overview of all the dependencies of the final results based on the implementation itself. Note that `leanblueprint` would not suffice, as it creates the dependency graph based on the user-specified dependencies in the latex blueprint and not on the Lean implementation itself.

Our next suggestion concerns *delaboration* in Lean, which is the inverse of *elaboration*. The process of elaboration takes user-facing syntax, which may be ambiguous in terms of, for example, the used notation, inferred types, and implicit arguments, and transforms it into Lean's core type theory, producing a typed expression. For instance, the type of `π + 2` can be different depending on the meaning of $\pi$ (which could be the transcendental real number close to 3, or the prime counting function), 2 (a numeric literal integer, a real number, or even the set $\{\{\}, \{\{\}\}\}$ in von Neumann ordinals), and $+$ (integer addition, real addition, or even direct sum of groups). The elaborator assigns Lean constants and variables to these syntactic objects using the current context: if a larger expression is `Real.sin (π + 2)`, then we expect a real number, which allows us to infer the meaning of all the syntactic objects. In contrast to elaboration, delaboration transforms a parsed tree back into a syntactic representation and is used as part of the interactive development process, as opposed to program execution. It is not to be confused with printing a value, or `ToRepr` in Lean, which transforms data into a human-readable representation. A *delaborator* is then a code that represents a particular class of expressions in a more readable form.

While working on our project, we discovered several shortcomings of the current delaborator system.

First, `Subtype` has a delaborator producing `⟨i, hi⟩`, which is more readable than the default `(Subtype.mk i hi)`, especially when there are many occurrences within a complex expression. Unfortunately, the custom `Subtype` delaborator is incompatible with the `pp.structureInstances=false` setting, which makes other structures more readable. It would be nice to be able to simultaneously use explicit constructors for some types and specialized delaborators for other types.

Another limitation of the current delaborator is its handling of dot notation. For example, given a list `l`, its length `List.length l` is displayed in the Infoview using dot notation as `l.length`. However, the delaborator does not always take advantage of the dot notation and sometimes displays the longer version instead, in particular for `private` declarations[9] and for `Function.something` or `Subtype.something` declarations. The way we addressed this in our implementation was by implementing custom delaborators, such as

```
@[app_unexpander Matrix.toCanonicalSigning]
private def Matrix.toCanonicalSigning_unexpand : Lean.PrettyPrinter.Unexpander
  | `($_ $Q) => `($(Q).$(Lean.mkIdent `toCanonicalSigning))
  | _ => throw ()
```

It would be ideal if delaboration consistently took advantage of dot notation out of the box, as implementing unexpanders manually clutters the code and is prone to errors.

### A.2 Suggestions for Mathlib

We extensively utilized Mathlib in our project, especially its support for matrices, linear independence, and matroids, and we saw three areas for potential improvement. First, the support for block matrices could be extended. Second, we noticed that the API for constructing single-row and single-column matrices has been getting less convenient over time. During our work on the project, converting a vector `r` to a single-row matrix containing `r` evolved from `Matrix.row r`, to `Matrix.row Unit r`, to `Matrix.replicateRow Unit r`. To make this more concise in our implementation, we introduced custom notation for such operations, and we wish the syntax in Mathlib was simplified. Last but not least, we found that the `variable` command is sometimes exploited in the source code of Mathlib to hide explicit arguments. We would suggest against using this command this way, to make the source code of the library easier to read and use.

### A.3 Suggestions for Aesop

The tactic `aesop` is very powerful and allowed us to prove quite a few lemmas automatically. Unfortunately, in certain situations `aesop` does not close the goal even when it most likely should. For example, it never uses `Ne.symm`, so it could happen that `aesop` can close the goal only with $x \neq y$ provided in the context, but not with $y \neq x$. Another downside is that `aesop` does not make use of the `fin_cases` tactic. This could be extremely beneficial whenever there is just a handful of cases that need to be considered. For example, the workload of analyzing `Fin 3` values is similar to analyzing all possible cases of a term from an inductive type with three constructors, yet the former is not supported.

[9] We created a Request for Comments https://github.com/leanprover/lean4/pull/10122 and, a few days before this paper's submission, this issue has been solved.